\documentclass[11pt,twoside,leqno]{amsart}
\usepackage{amsmath,amscd,amssymb,amsfonts}
\usepackage[all]{xy}

\renewcommand{\O}{{\mathcal{O}}}

\newcommand{\N}{\mathbb{N}}
\newcommand{\Z}{\mathbb{Z}}

\newcommand{\C}{\mathbb{C}}
\newcommand{\G}{\mathbb{G}}

\renewcommand{\P}{\mathbb{P}}

\newcommand{\g}{\mathfrak{g}}

\newcommand{\eps}{\varepsilon}

\newcommand{\be}{\begin{equation*}}
\newcommand{\ee}{\end{equation*}}
\newcommand{\ben}{\begin{equation}}
\newcommand{\een}{\end{equation}}

\newcommand{\ch}{\mathop{\mathrm{ch}}}

\newcommand{\td}{\mathop{\mathrm{td}}}

\newcommand{\RepD}{\mathop{\underline{\smash{\mathrm{Rep}}}}}
\newcommand{\Rep}{\mathop{\mathrm{Rep}}}

\title{Dimensional interpolation 
and the Selberg integral}

\author{V. Golyshev, D. van Straten, and D. Zagier}

\begin{document}

\bibliographystyle{alpha}

\begin{abstract} We show that a version of dimensional interpolation for the Riemann--Roch--Hirzebruch formalism in the case of a grassmannian leads to an expression for the Euler characteristic of line bundles in terms of a Selberg integral. We propose a way to interpolate higher Bessel equations, their wedge powers, and monodromies thereof to non--integer orders, and link the result with the dimensional interpolation of the RRH formalism in the spirit of the gamma conjectures.    

\end{abstract}

\maketitle

\bigskip
The dimensions of spaces of 
sections of certain ample bundles on homogeneous 
spaces such as the grassmannian $G(k,N+k)$ of $k$--spaces
in  $\C^{N+k}$
can easily be interpolated as functions of the  variable $N$ since 
they depend on $N$ polynomially. 
Deligne interpolated the spaces themselves to objects of a certain tensor category $\RepD(GL_t)$, where 
$t$ should be thought of as $N+k$. For instance, 
beginning with the polynomial interpolation of
$\chi_{\P^{N}} ({\mathcal{O}} (n))$ as ${N+n \choose n }={t-1+n \choose n }$, 
one can take a step further and interpret
$H (\P^N,{\mathcal{O}} (n))$ as $\mathop{\mathrm{Sym}}^n V_t^*$ in 
Deligne's category. We will leave $\P ^{N} = G(1,N+1)  = \P ^{t-1}$ itself undefined, trying instead to operate with its vestiges in a consistent manner.

Different levels of interpolation appear naturally 
in this framework.  For $t$ a natural number, the dimensions of actual objects in the usual category of representations $\Rep(GL_t)$  are given by the Weyl character formula and
interpolate easily as functions of the highest weight. But how do the individual ingredients of Weyl's formula interpolate as functions of the length of the Dynkin diagram? Or, by Hirzebruch, the Euler characteristic of a vector bundle on a homogeneous 
space is the result of pairing up the Chern character 
with the Todd genus of the space; can both be interpolated 
naturally in such a way that the result of 
the pairing still behaves polynomially? 

A level deeper, the Riemann--Roch--Hirzebruch numerology of projective spaces can be linked to the 
monodromy (or its version adapted to irregular connections) 
of the higher 
Bessel equations 

\begin{equation*}
(D^{N+1} - z^{N+1}) \Psi (z) = 0, \quad N \in \N,
\end{equation*}     
where 
$D=z \frac{d}{dz}.$ 
Does the  interpolated `differential equation'
\begin{equation*}
(D^t - z^t) \Psi (z) = 0,
\end{equation*}
have rudiments of monodromy? Is this monodromy 
related to an interpolation of the RRH theorem? 

Preliminary results 
show that the answers to these questions are positive. 
This short note should be viewed as a mere announcement of a more detailed version: we barely indicate the direction we are going
by explaining a link between dimensional interpolation of the RRH formalism for line bundles on grassmannians and the Selberg integral (in~\bf 3\rm~below)
and showing how a version of the gamma conjecture for grassmannians 
might remain true (in \bf 4\rm). 	
  
\medskip
We remark that several attempts to interpolate dimension as a continuous 
variable have been made in the past.
Hausdorff defined fractional dimension $d$
in terms of the asymptotics of the number
$\mathcal{N}(R)$ of balls of radius $R$ needed to cover $X$:
\[ \mathcal{N}(R) \sim \left(\frac{1}{R}\right)^d.\]
J. von Neumann defined fractional dimensions in his {\em Continuous Geometry} 
in terms of projectors in operator algebras. In quantum field theory, 
dimensional regularization was introduced
by t'Hooft and Veltman \cite{HV1972} and independently Bollini and 
Giambiagi \cite{BG1972}. The main observation is that the Feynman 
integral that belongs to a given Feynman diagram 
`can be taken in any dimension $D$', basically because its integrand
only contains (after Wick rotation) scalar products in euclidean space, 
\cite{Etingof1999}.
For a review of various approaches to dimension, and specifically, a link with 
the theory of modular forms, we refer the reader to \cite{Man2006}. 

Deligne's category, which we introduce in the next section, is not explicitly used anywhere in the paper. Nevertheless, its relevance should be clear from the context.

\subsection{Deligne's category  $\RepD(GL_t)$}
From the perspective of the
Killing-Cartan-Weyl classification of simple Lie algebras
and their representation theory in terms of highest weights, root systems, Weyl groups and associated combinatorics
it is not so easy to 
understand the extreme uniformity in the representation theory 
that exists among different Lie groups. With possible application 
to a universal Chern-Simons type knot invariant in mind, P. Vogel \cite{Vog1999}
tried to define a universal Lie algebra, $\g(\alpha:\beta:\gamma)$ 
depending on three  {\em Vogel parameters} that determine a point 
$(\alpha:\beta:\gamma)$ in the {\em Vogel plane}, in which all simple 
Lie algebras find their place. The dimension of the Lie algebra 
$\g(\alpha:\beta:\gamma)$ is given by a universal rational expression
\be
\dim \g(\alpha:\beta:\gamma)\, = \, \frac{(\alpha-2t)(\beta-2t)(\gamma-2t)}{\alpha\beta\gamma},\qquad
t=\alpha+\beta+\gamma ,
\ee and similar universal rational formulas can be
given for the dimensions of irreducible constituents of $S^2\g, S^3\g$ and 
$S^4\g$. Although the current status of Vogel's suggestions is unclear to us, 
these ideas have led to many interesting developments, such as the discovery of 
$E_{7\frac{1}{2}}$ by Landsberg and Manivel, 
\cite{LM2002}, \cite{LM2004}, \cite{LM2006}, \cite{LM2006a},    \cite{LM2006a}.

In order to interpolate within the classical $A,B,C,D$ series of Lie algebras, 
Deligne has defined $\otimes$-categories
\[ \RepD(GL_t),\;\; \RepD(O_t), \]
where $t$ is a parameter that can take on any complex value. 
(The category $\RepD(Sp_{2t})$ 
 is usually not discussed as it can be 
expressed easily in terms of the category $\RepD(O_T)$ with $T=-2t$.)
If $n$ is an
integer, there are natural surjective functors
\[\RepD(GL_n) \to \Rep(GL_n)\]
In the tannakian setup one would attempt to reconstruct a group $G$ from its
$\otimes$-category of representations $\Rep(G)$ using a fibre functor to
the $\otimes$-category $Vect$ of vector spaces, but Deligne's category has no fibre functor and is not tannakian, or, in general, even abelian. (However, when $t$ is not an integer, the category \emph{is}
abelian semisimple.)  

According to the axioms, in an arbitrary rigid $\otimes$-category $\mathcal{R}$ there exist a unit object~${\bf 1}$ and 
canonical evaluation and coevaluation morphisms
\[ \epsilon: V \otimes V^* \to {\bf 1},\qquad \delta: {\bf 1} \to V \otimes V^*\] 
so that we can assign to any object a dimension by setting
\[ \dim V =\epsilon \circ \delta  \in \mathop{\mathrm{End}}({\bf 1}) \in \C. \]

A simple diagrammatic description of $\RepD(GL_t)$ can be found in
\cite{CW2012}. One first constructs a skeletal category ${\RepD\,}_0(GL_t)$, 
whose objects are words in the alphabet $\{\bullet, \circ\}$. The letter
$\bullet$ corresponds to the fundamental representation $V$ of $GL_t$, 
$\circ$ to its dual $V^*$. A~$\otimes$-structure is induced by concatenation 
of words. The space of morphisms between two such words is the $\C$-span of 
a set of admissible graphs, with vertices the circles and dots of the two 
words.  Such an admissible graph consists of edges between the letters of 
the two words. Each letter is contained in one edge. Such an edge connects 
different letters of the same word or the same letter if the words are 
different.
$$\vcenter
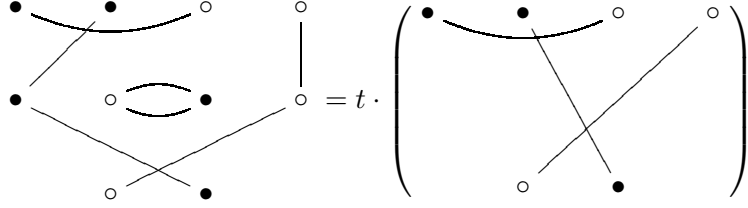
{ \xymatrix{
\bullet \ar@/_2ex/@{-}[rr] & \bullet \ar@{-}[ld] & \circ & \circ \ar@{-}[d] \\
\bullet \ar@{-}[rrd] & \circ  \ar@/^/@{-}[r] \ar@/_/@{-}[r]& \bullet & \circ \ar@{-}[lld] \\
		& \circ & \bullet & \\  
} } 
= t \cdot
\left( \vcenter{ \xymatrix{\bullet \ar@/_2ex/@{-}[rr] & \bullet \ar@{-}[ddr] & \circ & \circ \ar@{-}[ddll] \\
 &  &  \\
		& \circ & \bullet & \\  
}} \right) $$
The composition is juxtaposition of the two graphs, followed by
the elimination of loops, which results in a factor $t$.\\

Deligne's
category is now obtained by first forming its additive hull by introducing
formally direct sums  and then passing to the 
Karoubian hull, i.e. forming a category of pairs $(W,e)$, consisting of 
an object together with an idempotent:
\[\RepD (GL_t) =({\RepD\,}_0(GL_t)^{\text{add}})^\text{Karoubi}. \]

\bf Example. \rm  Consider the word $\bullet \bullet$ and 
the morphisms $\mathrm{Id}$ and $\mathrm{Swap}$ with the obvious meaning. 
One then can put
\[ S^2V=(\bullet \bullet, s), \;\;\wedge^2 V=(\bullet\bullet, a),\]
where
\[ s=\frac{1}{2}(\mathrm{Id}+\mathrm{Swap}),\;\;a=\frac{1}{2}(\mathrm{Id}-\mathrm{Swap})\]
so that in $\RepD(Gl_t)$ one has:
\[ V \otimes V=(\bullet \bullet ,Id)=S^2V \oplus \wedge^2V,\]
which upon taking dimensions is the identity
\[  t^2 = \frac{t(t+1)}{2}+\frac{t(t-1)}{2} .\]

\subsection{`Spaces of sections' as objects in Deligne's category and the beta integral.}
As above, we assume that $n$ is a natural number. Write $t=N+1$ and let $V_t=V$ be the fundamental object of $\RepD(GL_t)$ so that
$\dim V_t=t$. We do not define the projective space $\P =\P^N$, 
but we can pretend that, in the sense of Deligne, the space of global sections is
\[ H(\mathcal{O}_{\P}(n)) :=\mathop{\mathrm{Sym}^n}(V_t^*) \in \RepD(GL_t) .\]
Its dimension is then, as expected
\begin{equation}
\chi(\mathcal{O}_{\P}(n)) :=\dim H (\mathcal{O}_{\P}(n))={N+n \choose n} \label{chi-interpret},
\end{equation}
(interpreted in the obvious way as a polynomial
in $N$ if $N\not\in\Z$),
so that e.g.
\[\chi(\mathcal{O}_{\P^{1/2}}(2))=\frac{3}{8}.\]
The Poincar\'e series 
\[P(y):=\sum_{n=0}^{\infty} \chi(\mathcal{O}_{\P}(n)) y^n =\frac{1}{(1-y)^{N+1}},\]
which is consistent with the idea that $\dim V_t = N+1$.

\medskip
Returning to the question posed at the beginning, 
`is there a way to extend the interpolation of $\chi$
individually to the Chern and the Todd ingredients?', we reason as
follows. If $X$ is a smooth projective $n$-dimensional variety, and $E$ a vector 
bundle on $X$, then the Euler characteristic
\[\chi(X,E):=\sum_{i=0}^n (-1)^i\dim H^i(X,E)\]
can be expressed in terms of characteristic numbers
\[\chi(X,E)=\int_X \ch(E) \cdot \td(X) . \] 
Here the integral in the right hand side is usually
 interpreted 
as resulting from evaluating the cap product with the fundamental class $[X]$
on the cohomology algebra $H^*(X)$, and the Chern character and Todd class
are defined in terms of the Chern roots $x_i$ of $E$ and $y_i$ of $TX$:
\[\ch(E)=\sum_{i=1}^r e^{x_i}\,, \qquad \td(X)=\prod_{i=1}^n \frac{y_i}{1-e^{-y_i}} .\]
The cohomology ring of an $n$-dimensional projective space is
a truncated polynomial ring:
\[H^*(\P^N)=\Z[\xi]/(\xi^{N+1})\,, \qquad\xi=c_1(\mathcal{O}(1)),\]
and it is not directly clear how to make sense of this if $N$ is not an 
integer. Our tactic will be  to drop the relation 
\[\xi^{N+1}=0\] 
altogether, thinking instead of $\Z[\xi]$  as a Verma module over the  $sl_2$
of the Lefschetz theory, and replacing taking the  cap product 
with integration. 
As we will be integrating meromorphic functions in $\xi$,
the polynomial ring is too small, and we put  
\[ \hat{H}(\P) :=\Z[[s]] \supset \Z[s] .\] 
One has
\[ \chi(\mathcal{O}(n))=e^{n\xi}\,, \qquad \td(\P)=\left(\frac{\xi}{1-e^{-\xi}}\right)^{N+1}, \]
so Hirzebruch-Riemann-Roch reads
\[\chi(\mathcal{O}(n))=\left[e^{n\xi} \left(\frac{\xi}{1-e^{-\xi}}\right)^{N+1}\right]_N\]
where $[...]_N$ is the coefficient at $\xi^N$ in a series. 
This can be expressed analytically as a residue integral
along a small circle around the origin:
\be
\chi(\mathcal{O}(n))=\frac{1}{2\pi i}\oint e^{n \xi}\left(\frac{\xi}{1-e^{-\xi}}\right)^{N+1}\frac{d\xi}{\xi^{N+1}} .
\ee
As it stands, it cannot be extended to non-integer
$N$ since the factor $(1- e^{-\xi})^{-N-1}$ is not univalued on the circle.  The usual way to adapt it is to  consider, for $n \ge 0$, the integral along the   path 
going from $- \infty - i \eps$ to $ - i \eps$, making a half--turn round the origin and going back, and choosing the standard branch of the logarithm.   Because of the change in the argument this integral is equal to  
\be
J(N,n) = 
\frac{e^{2 \pi i (N+1)}-1}{2 \pi i} 
\int_{-\infty}^0 \frac{e^{n \xi}}{(1-e^{-\xi})^{N+1}} 
d \xi , 
\ee
or, after the substitution $s=e^{\xi}$,
\be
J(N,n) = 
\frac{e^{2 \pi i (N+1)}-1}{2 \pi i} 
\int_0^1 s^{n-1} (1-1/s)^{-N-1} ds
=
\frac{\sin \pi (N+1)}{ \pi } 
\int_0^1 s^{n+N} (1-s)^{-N-1} ds .
\ee
Using Euler's formulas
\ben
\Gamma(x)\Gamma(1-x) =\frac{\pi}{\sin \pi x} \,,
\label{gamma-one-minus-argument}
\een
\begin{equation}
\int_0^1 s^{\alpha-1}(1-s)^{\beta-1} ds = \frac{\Gamma(\alpha)\Gamma(\beta)}{\Gamma(\alpha+\beta)} \,, 
\label{beta-integral}
\end{equation}
and
\be
\frac{\Gamma(N+n+1)}{\Gamma(n+1) \Gamma(N+1)}
=
{N+n \choose n} \,,
\ee
we arrive at a version of RRH `with integrals':

\medskip

\noindent \bf Proposition 1. \rm Let $n \in \N$. Assume $\mathop{\mathrm{Re}} N  < 0, \, N \notin \Z$. Interpret the Euler characteristic of~$\P ^N$ via formula \eqref{chi-interpret}. Then  
\be \label{little-propo}
\chi_\P(\O (n)) = \frac{e^{2 \pi i (N+1)}-1}{2 \pi i} 
\int_{-\infty}^0 \frac{e^{n \xi}}{(1-e^{-\xi})^{N+1}} 
d \xi .
\ee
\qed

\bigskip

%

\medskip

\subsection{The grassmannian and the Selberg integral.} For $\P^N$, we ended up with the beta function, a one-dimensional integral, as the cohomology ring 
is generated by a single class~$\xi$. In the cases  where the cohomology ring is generated
by $k$ elements, for example the grassmannian $G(k,N+k)$,
we would like to see a $k$-dimensional integral appear in a natural way.
For $N \in \N$ the cohomology ring of the grassmannian $\G:=G(k,N+k)$ 
is given by
\[H^*(G(k,N+k))=\C[s_1,s_2,\ldots,s_k]/(q_{N+1},q_{N+2}, \dots,     q_{N+k}),\]
where the $s_i$ are the Chern classes of the universal rank $k$ sub-bundle
and $q_i=c_i(Q)$ are formally the Chern classes of the universal quotient bundle $Q$ (so that the generating series of $q$'s is inverse to that of $s$'s). 
In the same vein as before, we set:
\ben
\hat{H}^*(\G):=\C[[s_1,s_2,\ldots,s_k]]=\C[[x_1,x_2,\ldots,x_k]]^{S_k} \label{drop-rel}
\een
by dropping the relations. A $\C$-basis of this ring
given by the Schur polynomials
\[\sigma_{\lambda} :=\frac{\det(x_i^{\lambda_j+k-j})}{\det(x_i^{k-j})}\]		
where $\lambda=(\lambda_1,\lambda_2,\ldots,\lambda_k)$ 
is an arbitrary Young diagram with at most $k$ rows.
There is a Satake--type map for the extended cohomology: 
\[ \mathrm{Sat}: \hat{H}(\G) \to \wedge^k \hat{H}(\P) \]
obtained from the Young diagram by `wedging its rows':
\[\sigma_{\lambda} \mapsto  s^{\lambda_1+k-1} \wedge s^{\lambda_2+k-2}\wedge \ldots \wedge s^{\lambda_k-1}. \] 
We are therefore seeking an expression for the values of the Hilbert polynomial of $G(k,N)$ in terms of a $k$--dimensional integral of the beta type involving $k$--wedging.

Euler's beta integral
\eqref{beta-integral}
has several generalizations.  Selberg introduced \cite{Selberg1944}
an integral \cite{FW2008} over the $k$-dimensional cube
\begin{equation*}
S(\alpha, \beta,\gamma, k):=\int_0^1\ldots\int_0^1 (s_1 s_2\ldots s_k)^{\alpha-1}((1-s_1)(1-s_2)\ldots(1-s_k))^{\beta-1}\Delta(s)^{2\gamma} ds_1ds_2\ldots ds_k 
\end{equation*} 
where
\[ \Delta(s)=\Delta(s_1,s_2,\ldots,s_k)=\prod_{i <j} (s_i-s_j) ,\]
and showed that it admits meromorphic continuation, which
we will also denote by $S$. 

\medskip
\noindent {\bf Proposition 2.}  For $k \in \mathbb{N},\,  n \in \mathbb{Z}_+$, let $\chi(\mathcal{O}_{\G}(n))$
denote the result of interpolating  the polynomial function $\chi(\mathcal{O}_{G (k,k+N)}(n))$  of the argument 
$N \in \N$ to $\C$.  
  One has
\begin{equation*}
\chi(\mathcal{O}_{\G}(n))= \frac{(-1)^{k(k-1)/2}}{k!}\left( \frac{\sin \pi(N+1)}{\pi}\right)^k S(n+N+1,-N-k+1;1,k) .
\end{equation*} \rm
{\sc Proof.} The shortest (but not the most transparent) way to see this is to use the expressions for the LHS and the RHS in terms of the product of gamma factors found by Littlewood and Selberg respectively. By Selberg, 
\begin{equation}
 S(\alpha,\beta,\gamma,k)=\prod_{i=0}^{k-1} \frac{\Gamma(\alpha+i\gamma)\Gamma(\beta+i\gamma)\Gamma(1+(i+1)\gamma)}
{\Gamma(\alpha+\beta+(k+i-1)\gamma) \Gamma(1+\gamma)} \label{Selberg-formula}.
\end{equation}
By Littlewood \cite{Lit1942}, for $N \in \Z_{>0}$ one has
\begin{equation*}
 \chi(\mathcal{O}_{G(k,k+N)}(n)) =\frac{{N+n \choose n} {N+n+1 \choose n+1} \ldots {N+n+(k-1) \choose n+(k-1)}}{
{N \choose 0} {N+1 \choose 1} \ldots {N+(k-1) \choose (k-1)}}, 
\end{equation*}
where there are $k$ factors at the top and the bottom. 
Rearranging the terms in the RHS of \eqref{Selberg-formula} and using 
\eqref{gamma-one-minus-argument},
we bring the $\Gamma$-factors that involve $\beta$ to the denominator in order to form the binomial coefficients at the expense of the sine factor.
 \qed

\bigskip
As an example, for $k=2$ and 
$N=-1/2$, we get
the Hilbert series 
\[ \sum_{k=0}^\infty \chi(\mathcal{O}_{G(2,3/2)}(n)) \, y^n = 1+6\, \frac{t}{16} +60\left(\frac{t}{16}\right)^2 + 700 \left(\frac{t}{16}
\right)^3+8820 \left(\frac{t}{16}
\right)^4  +\ldots
\]
which is no longer algebraic, but can be expressed in terms of elliptic functions.

More generally, one can consider a Selberg--type integral 
with an arbitrary symmetric function rather than the discriminant in the numerator and use separation of variables together with the Jacobi--Trudi formula in order to obtain similar expressions in terms of the gamma function in order to interpolate between the Euler characteristics of more general vector bundle on grassmannians (or the dimensions of highest weight
representations of $GL_{N+k}$).

\subsection{Towards a gamma conjecture in 
non--integral dimensions.} \label{gamma-phenomena} The by now standard 
predictions of mirror symmetry relate the RRH formalism 
on a Fano variety $F$ to the monodromy of its regularized 
quantum differential equation. It is expected that this
differential equation arises from the Gauss--Manin
connection  
in the middle cohomology of level hypersurfaces of a 
regular function $f$ defined on some quasiprojective 
variety (typically a Laurent 
polynomial on $\mathbb{G}_\mathrm{m}^{\, d}$), 
called in this case a Landau--Ginzburg model of 
$F$. By stationary phase, 
the monodromy of the Gauss--Manin connection in a pencil 
translates
into the asymptotic behavior of oscillatory 
integrals of the generic form
$I (z) = \int \exp (izf)\, d\mu (\mathbb{G}_\mathrm{m}^{\, d})$,
which satisfy the quantum differential equation of $F$, 
this time without the word `regularized'. 
The asymptotics 
are given by Laplace integrals computed at the critical points,
and the critical values of $f$ are the exponents occurring
in the oscillatory integrals $I_i(z)$ that have `pure' 
asymptotic behavior in sectors. 
One wants to express these pure asymptotics in terms of the Frobenius 
basis of solutions $\{ \Psi_i (z) \}$ around $z = 0$.
The gamma conjecture \cite{GGI2016} 
predicts that such an expression 
for the highest--growth asymptotic (arising from the critical
value next to infinity) will give the `gamma--half' of the Todd 
genus and therefore effectively encode the Hilbert
polynomial of $F$ with respect to the anticanonical bundle. 
At first sight, none of this seems capable of surviving in non--integer dimensions.
Yet, to return to the example of $G(2,N+2)$, define the numbers $c_j$ and $d_j$ by the expansions 
\be
\Gamma_\P^{(0)} (\eps) = \Gamma (1+\eps)^{N+2} = \sum_{j=0}^\infty d_j \eps^j ,
\ee
\be
\Gamma_\P^{(1)} (\eps) = \Gamma (1+\eps)^{N+2} e^{2 \pi i \eps} = \sum_{j=0}^\infty c_j \eps^j.
\ee
Put
\be
F(\eps,z) = \sum_{k=0}^{\infty} \frac{z^{l+\eps}}{\Gamma(1+l+\eps)^{N+2}}
\ee
and
\be
\Psi (\eps,z) = \Gamma_\P (\eps) F (\eps, z) = \sum_{k=0}^\infty \Psi_k (z) \eps^k.
\ee

\medskip
\noindent \bf Claim \rm (rudimentary gamma conjecture). For fixed $N > 2$ and $i, \, j$  in a box of at least some moderate size,  
one should have
\be
\label{claim-grass}
\lim_{z \to - \infty} \frac{\Psi_i (z) \Psi'_j (z) - \Psi_j (z) \Psi'_i (z)}{\Psi_1 (z) \Psi'_0 (z) - \Psi_0 (z) \Psi'_1 (z)}
= \frac{c_i d_j - c_j d_i}{c_1 d_0 - c_0 d_1} .
\ee

\bigskip
\bigskip

\noindent The LHS and RHS mimic, in the setup of formula \eqref{drop-rel}, the $\sigma_{[j-1,i]}$-coefficients in the expansion of the `principal asymptotic class' and the gamma class of the usual grassmannian: in the case when $N \in \N$ and $0 \le  i,j  \le N$ one would use the identification of $2$--Wronskians of a fundamental matrix of solutions to a higher Bessel equation with homology classes of $G(2,N+2)$.  Preliminary considerations together with numerical evidence suggest that the claim has a good chance to be true, as well as its versions for $G(k,N+k)$ with $k > 2$.

\bigskip

%

\bigskip
\bigskip

The first--named author is grateful to Yuri Manin and Vasily Pestun for stimulating discussions. We thank Hartmut Monien for pointing us to \cite{FW2008}.

\bigskip

\nocite{MV2017}
\nocite{GM2014}
\nocite{Bra2013}
\nocite{BS2013}
\nocite{Etingof1999}
\nocite{Etingof2014}
\nocite{Etingof2016}
\nocite{EGNO2015}
\nocite{FW2008}
\nocite{Man2006}
\nocite{Opd1999}
\nocite{Man1985}
\nocite{Lit1942}
\nocite{Lit1943}
\nocite{BD2016}
\nocite{LM2002}
\nocite{LM2004}
\nocite{LM2006}
\nocite{LM2006a}

\nocite{GW2011}
\nocite{Del2002}
\nocite{}
\nocite{}
\nocite{}
\nocite{}
\nocite{}
\nocite{}

\vskip1.3cm \noindent
{\sc Algebra and Number Theory Laboratory \newline\noindent
Institute for Information Transmission Problems \newline\noindent
Bolshoi Karetny 19, Moscow 127994, Russia}

\smallskip\noindent
golyshev@mccme.ru

\bigskip   

\noindent
{\sc Institut f\"ur Mathematik \newline
Johannes Gutenberg-Universit\"at \newline
Staudingerweg 9, 4. OG
55128 Mainz, Germany
}

\smallskip\noindent
straten@mathematik.uni-mainz.de

\bigskip   

\noindent
{\sc Max Planck Institut f\"ur Mathematik\newline\noindent
Vivatsgasse 7, 53111 Bonn, Germany}\newline\noindent
and \newline\noindent
{\sc International Centre for Theoretical Physics\newline\noindent
Via Miramare, Trieste, Italy}

\smallskip\noindent
dbz@mpim-bonn.mpg.de

\end{document}